\documentclass[11pt,english]{article}

\usepackage{graphicx,amssymb,amsmath,empheq,bbm,float,subfigure,nicefrac,bm,setspace,parskip,enumerate,authblk,xcolor,booktabs}
\usepackage[numbers, sort&compress]{natbib}
\definecolor{midnightblue}{cmyk}{1,1,0,0.1}
\definecolor{forestgreen}{cmyk}{0.76,0,0.26,0.5}

\usepackage[a4paper]{geometry}
\usepackage{hyperref}
\hypersetup{
    colorlinks=true,       
    linkcolor=midnightblue,          
    citecolor=magenta,        
    filecolor=midnightblue,      
    urlcolor=midnightblue,          
}

\providecommand{\R}{{\ensuremath{\mathbb{R}}}}
\providecommand{\E}{{\ensuremath{\mathbbm{E}}}}
\renewcommand{\P}{\mathbbm{P}}
\DeclareMathOperator{\Tr}{Tr}

\makeatletter
\usepackage{babel}
\usepackage{csquotes}

\title{Solving High-Dimensional Partial Differential Equations Using Deep Learning}

\author[1]{Jiequn Han}
\author[2]{Arnulf Jentzen}
\author[4,3,1]{Weinan E\thanks{weinan@math.princeton.edu}}
\affil[1]{Program in Applied and Computational Mathematics, \protect\\ Princeton University, Princeton, NJ 08544, USA}
\affil[2]{Department of Mathematics, ETH Z\"urich, R\"amistrasse 101, 8092 Z\"urich,
Switzerland}
\affil[3]{Department of Mathematics, Princeton University, Princeton, NJ 08544, USA}
\affil[4]{Beijing Institute of Big Data Research,
Beijing, 100871, China}
\begin{document}
\date{\vspace{-1ex}}
\maketitle

\begin{abstract}
Developing algorithms for solving high-dimensional partial differential equations (PDEs) has been 
an exceedingly difficult task for a long time, 
due to the notoriously difficult problem known as the
 ``curse of dimensionality''. 
This paper introduces a deep learning-based approach that can handle general
high-dimensional parabolic PDEs. To this end, the PDEs are reformulated
using backward stochastic differential equations and the gradient of the unknown solution is approximated by neural networks, very much in the spirit of deep reinforcement learning
with the gradient acting as the policy function. 
Numerical results on examples including the nonlinear Black-Scholes equation, 
the Hamilton-Jacobi-Bellman equation, and the Allen-Cahn equation 
suggest that the proposed algorithm is quite effective 
in high dimensions, in terms of both accuracy and cost. 
This opens up new possibilities in economics, finance, operational research, and physics,
by considering all participating agents, assets, resources, or particles together
at the same time, instead of making ad hoc assumptions on their inter-relationships.
\end{abstract}

\section{\label{sec1}Introduction}
Partial differential equations (PDEs) are among the most ubiquitous tools 
used in modeling problems in nature.
Some of the most important ones are naturally formulated as PDEs in high dimensions. 
Well-known examples include:
\begin{enumerate}
  \item The Schr\"odinger equation in quantum many-body problem. In this case the 
dimensionality of the PDE is roughly three times the number of electrons or quantum particles
in the system.

  \item The nonlinear Black-Scholes equation for pricing financial derivatives, 
in which the dimensionality of the PDE is the number of underlying financial assets 
under consideration.

  \item The Hamilton-Jacobi-Bellman equation in dynamic programming. 
  In a game theory setting with multiple agents, the dimensionality goes up linearly with 
  the number of agents.  Similarly, in a resource allocation problem, 
  the dimensionality goes up linearly with the number of devices and resources.
\end{enumerate}

As elegant as these PDE models are, their practical use has proven to be 
very limited
due to the curse of dimensionality \cite{Bellman1957}: the computational cost for solving
them goes up exponentially with the dimensionality.

Another area where the curse of dimensionality has been an essential obstacle is
machine learning and data analysis, where the complexity of nonlinear regression models,
for example, goes up exponentially with the dimensionality.
In both cases the essential problem we face
is how to represent or approximate a nonlinear function in high dimensions.
The traditional approach, by building functions using polynomials, piecewise polynomials,
wavelets, or other basis functions, is bound to run into the curse of dimensionality 
problem.

In recent years a new class of techniques, the deep neural network model, has shown remarkable success in artificial intelligence
(see, e.g., \cite{Goodfellow2016,LeCun2015,Krizhevsky2012,Hinton2012,Silver2016}).
 Neural network is an old idea but recent
experience has shown that deep networks with many layers seem to do a
surprisingly good job in modeling complicated data sets.
In terms of representing functions, the neural network model is compositional:
it uses compositions of simple functions to approximate
complicated ones.  In contrast, the approach of classical approximation theory is usually additive.
Mathematically, there are universal approximation theorems stating that a single hidden layer neural network can approximate a wide class of functions on compact subsets
(see, e.g., survey \cite{Pinkus1999} and the references therein), even though we still lack a theoretical framework for explaining the seemingly unreasonable effectiveness of multilayer neural networks, which are widely employed nowadays. Despite this, the practical success of deep neural networks in artificial intelligence has been very astonishing and encourages applications to
other problems where the curse of dimensionality has been a tormenting issue.

In this paper, we extend the power of deep neural networks to another dimension
by developing a strategy for solving a large class of high-dimensional
nonlinear PDEs using deep learning. The class of PDEs that we deal
with are (nonlinear) parabolic PDEs.  
Special cases include the Black-Scholes equation and the Hamilton-Jacobi-Bellman equation.
To do so, we make use of the reformulation of these PDEs as backward stochastic
differential equations (BSDEs) (see, e.g., \cite{Pardoux1992,Pardoux1999}) and 
approximate the gradient of the solution using deep neural networks. 
The methodology bears some resemblance to deep reinforcement learning with the BSDE playing
the role of model-based reinforcement learning (or control theory models)
and the gradient of the solution
playing the role of policy function.
Numerical examples manifest that the proposed algorithm is quite satisfactory in both
accuracy and computational cost.

Due to the  ``curse of dimensionality'', there are only a very limited number of cases where practical high-dimensional algorithms have been developed in literature.
For linear parabolic PDEs, one can use the Feynman-Kac formula and Monte Carlo methods to develop efficient algorithms
to evaluate solutions at any given space-time locations.  
For a class of inviscid Hamilton-Jacobi equations, 
Darbon \& Osher have recently developed an effective algorithm 
in the high-dimensional case (see \cite{Darbon2016}), based on the Hopf formula for the Hamilton-Jacobi equations.
A general algorithm for nonlinear parabolic PDEs
based on the multilevel decomposition of
Picard iteration is developed in \cite{Weinan2016} and has been 
shown to be quite efficient on a number of examples in finance and physics.
The branching diffusion method has been proposed in \cite{Henry-Labordere2012,Henry-Labordere2014}, which exploits  the fact that  solutions of semilinear PDEs with polynomial nonlinearity can be represented as an expectation of a functional of  branching diffusion processes. This method does not suffer from the curse of dimensionality, but still has limited applicability  due to the blow up of approximated solutions in finite time.

The starting point of the present paper is deep learning.  It should be stressed that even though  deep learning has been a very successful tool for a number  of applications, adapting it to the current setting with practical success is still a highly non-trivial task.  Here by using the reformulation of BSDEs, we are able to cast the problem of solving PDEs as a learning problem and we design a deep learning framework that fits naturally to that setting. This has proven to be quite successful in practice. 


\section*{\label{sec2}Methodology}
We consider a general class of PDEs known as semilinear parabolic PDEs.
These PDEs can be represented as follows:
\begin{equation}
  \label{eq:PDE}
  \frac{ \partial u}{ \partial t } ( t, x )
  + \frac{1}{2} \! \Tr\!\Big( \sigma\sigma^{\operatorname{T}}(t,x)(\mbox{Hess}_x u) ( t, x ) \Big)
  +\nabla u( t, x )\cdot \mu( t, x )
  +f\big( t, x, u(t,x), \sigma^{\operatorname{T}}( t, x ) \nabla u( t, x ) \big) = 0 
\end{equation}
with some specified terminal condition $u(T,x) = g(x)$. 
Here $t$ and $x$ represent the time 
and $d$-dimensional space variable respectively,
$ \mu $ is a known vector-valued function, 
$\sigma $ is a known $d \times d$ matrix-valued function, 
$\sigma^{\operatorname{T}}$ denotes the transpose associated to $\sigma$, 
$\nabla u$ and $\mbox{Hess}_x u$ denote the gradient and the Hessian of function $u$ respect to $x$,
$\Tr$ denotes the trace of a matrix, and $ f $ is a known nonlinear function.
To fix ideas, we are interested in the solution at $ t = 0 $, $ x = \xi $ 
for some vector $ \xi \in \R^d $.

Let $ \{ W_t \}_{ t \in [0,T] } $ be a $ d$-dimensional Brownian motion and
$ \{ X_t \}_{ t \in [0,T] } $ be a $d$-dimensional stochastic process 
which satisfies
\begin{equation}
  \label{eq:FSDE}
  X_t = \xi + \int_0^t\mu( s, X_s )\, ds +\int_0^t \sigma( s, X_s ) \, dW_s.
\end{equation}
Then the solution of \eqref{eq:PDE} satisfies the following BSDE 
(cf., e.g., \cite{Pardoux1992,Pardoux1999}):
\begin{equation}
  \label{eq:BSDE_explicit}
  u(t, X_t) = 
  u(0, X_0) 
  - 
  \int_0^t f\big( 
    s, X_s, u(s,X_s), \sigma^{ \operatorname{T} }( s, X_s ) \, \nabla u( s, X_s ) 
  \big) \, ds
  + \int_0^t [ \nabla u( s, X_s ) ]^{ \operatorname{T} } \,\sigma( s, X_s )\, dW_s.
\end{equation}
We refer to the Section Materials and Methods for further explanation of \eqref{eq:BSDE_explicit}.

To derive a numerical algorithm to compute $u(0, X_0)$, 
we treat $u(0,X_{0}) \approx \theta_{u_0}, \nabla u(0,X_{0})\approx\theta_{\nabla u_0}$ 
as parameters in the model and 
view \eqref{eq:BSDE_explicit} as a way of computing the values of $u$ at the terminal 
time $T$, knowing $u(0, X_0)$ and $\nabla u(t, X_t)$.
We apply a temporal discretization to \eqref{eq:FSDE}--\eqref{eq:BSDE_explicit}. 
Given a partition of the time interval $[0, T]$:
$0 = t_0 < t_1 < \ldots < t_N = T$, we consider the simple Euler scheme for $n = 1, \dots, N-1$:
\begin{equation}
\label{eq:approx_Xt}
  X_{ t_{ n + 1 } } - X_{ t_n } \approx
  \mu( t_n, X_{ t_n } )\,\Delta t_n +
  \sigma( t_n, X_{ t_n } )\,\Delta W_n,
\end{equation}
and
\begin{equation}
\label{eq:approx_ut}
\begin{split}
  u( t_{ n + 1 },X_{t_{n+1}})
& \approx u( t_{ n },X_{t_{n}}) -
  f\big( 
    t_n, X_{ t_n }, u( t_{ n },X_{t_{n}}) , \sigma^{ \operatorname{T} }( t_n, X_{ t_n } )
    \,
    \nabla u( t_n, X_{ t_n } ) 
  \big)
  \,( t_{ n + 1 } - t_n )\\
&\quad +
    [ \nabla u( t_n, X_{ t_n } ) ]^{ \operatorname{T} } 
    \, \sigma( t_n, X_{ t_n } ) \,
    ( W_{ t_{ n + 1 } }-W_{ t_n }),
\end{split}
\end{equation}
where
\begin{equation}
  \label{eq:time_discret}
  \Delta t_n = t_{ n + 1 } - t_n, \quad \Delta W_n = W_{ t_{ n + 1 } }-W_{ t_n }.
\end{equation}
Given this temporal discretization, the path 
$
  \{ X_{ t_n } \}_{ 0 \leq n \leq N }
$
can be easily sampled using \eqref{eq:approx_Xt}. Our key step next is to approximate the function
$ x \mapsto \sigma^{ \operatorname{T} }(t,x) \, \nabla u(t,x)$
at each time step $ t = t_n $ 
by a multilayer feedforward neural network
\begin{equation}
\label{eq:approx_grad}
  \sigma^{ \operatorname{T} }( t_n , X_{ t_n } ) 
  \, 
  \nabla u( t_n, X_{ t_n } ) 
  =
  ( \sigma^{ \operatorname{T} } \nabla u )( t_n , X_{ t_n } )
  \approx ( \sigma^{ \operatorname{T} } \nabla u )( t_n, X_{ t_n } | \theta_n),
\end{equation}
for $n = 1, \dots, N-1$, where $\theta_n$ denotes parameters of the neural network approximating 
$ x \mapsto \sigma^{ \operatorname{T} }(t,x) \, \nabla u(t,x)$
at $ t = t_n $.

Thereafter, we stack all the sub-networks in \eqref{eq:approx_grad} together to form a deep neural 
network as a whole, based on the summation of \eqref{eq:approx_ut} over $n = 1, \dots, N-1$. 
Specifically, this network takes the paths 
$\{ X_{ t_n } \}_{ 0 \leq n \leq N }$ and 
$\{ W_{ t_n } \}_{ 0 \leq n \leq N }$ 
as the input data and gives the final output, 
denoted by 
$
  \hat{u}( 
    \{ { X_{ t_n } } \}_{ 0 \leq n \leq N } , 
    \{ W_{ t_n } \}_{ 0 \leq n \leq N } 
  ) 
$, 
as an approximation of 
$u( t_N, X_{ t_N } )$. 
We refer to the Section Materials and Methods for more details on the architecture of the neural network.
The difference in the matching of given terminal condition can be used to define the expected loss function 
\begin{equation}
  l(\theta) = 
  \E\Big[
    \big|
      g( X_{ t_N } ) -
      \hat{u}\big(
        \{ X_{ t_n } \}_{ 0 \leq n \leq N } ,
        \{ W_{ t_n } \}_{ 0 \leq n \leq N }
      \big)
    \big|^2
  \Big].
\end{equation}
The total set of parameters are:
$\theta=\{\theta_{u_0},\theta_{\nabla u_0},\theta_1,\dots,\theta_{N-1}\}$.

We can now use a stochastic gradient descent-type (SGD) algorithm to optimize the parameter $\theta$, just as in the standard training of deep neural networks.
In our numerical examples, we use the Adam optimizer \cite{Kingma2015}. 
See the Section Materials and Methods for more details on the training of the deep neural networks.
Since the BSDE is used as an essential tool, we call the
methodology introduced above \emph{deep BSDE method}.

\section*{\label{sec3}Examples}
\subsection*{Nonlinear Black-Scholes Equation with Default Risk}
A key issue in the trading of financial derivatives is to determine an appropriate fair price.
Black \& Scholes illustrated 
that the price $ u$ of a financial derivative satisfies a parabolic 
PDE, nowadays known as the Black-Scholes equation \cite{Black2012}.
The Black-Scholes model can be augmented to take into account several important
factors in real markets, including defaultable securities, 
higher interest rates for borrowing than for lending, 
transactions costs, uncertainties in the model parameters, etc. (see, e.g., \cite{Duffie1996,Bender2017,Bergman1995,Leland1985,Avellaneda1995}). 
Each of these effects results in a nonlinear contribution in the pricing model
(see, e.g., \cite{Bender2017,Crepey2013,Forsyth2001}). 
In particular, the credit crisis and the ongoing European sovereign debt crisis 
have highlighted the most basic risk that has been neglected in the original Black-Scholes model, 
the default risk \cite{Crepey2013}.

Ideally the pricing models should take into account
the whole basket of underlyings that the financial derivatives depend on,
resulting in high-dimensional nonlinear PDEs.
 However, existing pricing algorithms are unable to tackle
these problems generally due to the curse of dimensionality. 
To demonstrate the effectiveness of the deep BSDE method, we study
a special case of the recursive valuation model with default risk \cite{Duffie1996,Bender2017}. 
We consider the fair price of a European claim based on 100 underlying assets
 conditional on no default having occurred yet. 
When default of the claim's issuer occurs, the claim's holder only receives a fraction $\delta \in [0, 1)$ of the current value.
The possible default is modeled by the first jump time of a Poisson process with intensity $Q$, a decreasing function of the current value, i.e., the default becomes more likely when the claim's value is low. 
The value process can then be modeled by \eqref{eq:PDE} with the generator 
\begin{equation}
  f\big( t, x, u(t,x), \sigma^{\operatorname{T}}( t, x ) \nabla u( t, x ) \big)
  = 
  - \left( 1 - \delta \right) Q( u(t,x) ) \, u(t,x) 
  - R \, u(t,x)
\end{equation}
(see \cite{Duffie1996}), where $R$ is the interest rate of the risk-free asset. 
We assume that the underlying asset price moves as a geometric Brownian motion and choose 
the intensity function $Q$ as a piecewise-linear function of the current value with 
three regions ($v^h<v^l,\,\gamma^h>\gamma^l$):
\begin{equation}
  Q(y) = 
    \mathbbm{1}_{ ( - \infty, v^h) }(y)\,\gamma^h
    +
    \mathbbm{1}_{[ v^l, \infty )}( y )\,\gamma^l
    +
    \mathbbm{1}_{[ v^h, v^l )}(y )
    \left[
      \tfrac{ ( \gamma^h - \gamma^l ) }{ ( v^h - v^l ) }\left( y - v^h \right)+\gamma^h
    \right]
\end{equation}
(see \cite{Bender2017}).
The associated nonlinear Black-Scholes equation in $[0,\,T]\times\R^{100}$ becomes
\begin{equation}
  \label{eq:PDE_defaultrisk}
   \frac{ \partial u}{ \partial t } ( t, x )
  +
 \bar{\mu}x\cdot \nabla u(t,x)
  +
  \frac{ \bar{\sigma}^2 }{ 2 }\sum_{ i = 1 }^d| x_i |^2\,
\frac{ \partial^2 u}{ \partial x^2_i } (t,x)\\
- \left( 1 - \delta \right) Q( u(t,x) ) \, u(t,x) - R \, u(t,x)= 0.
\end{equation}

\begin{figure}[ht]
\centering
\includegraphics[width =8cm]{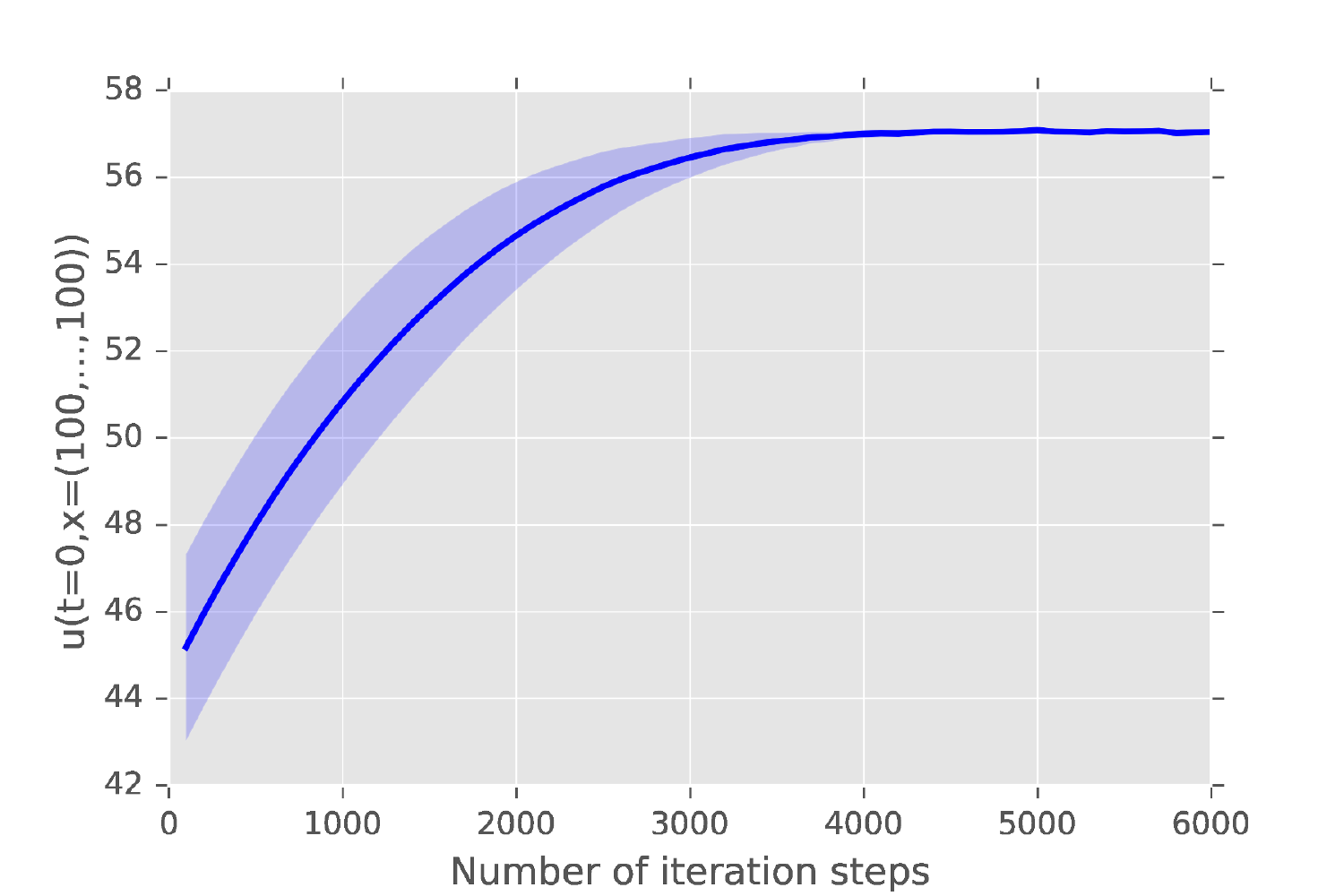}
\caption{Plot of $\theta_{u_0}$ as an approximation
of $u(t{=}0,x{=}(100,\dots,100))$
against the number of iteration steps
in the case of the $ 100 $-dimensional nonlinear Black-Scholes equation \eqref{eq:PDE_defaultrisk} 
with $ 40 $ equidistant time steps ($ N{=}40 $) and learning rate $ 0.008 $. 
The shaded area depicts the mean $\pm$ the standard deviation of $\theta_{u_0}$ as an approximation of $u(t{=}0,x{=}(100,\dots,100))$ for 5 independent runs.
The deep BSDE method achieves a relative error of size $ 0.46\% $ in a runtime of $ 1607 $ seconds.
\label{fig:defaultrisk}}
\end{figure}

We choose $T=1, \delta = 2/3,\,R=0.02,\,\bar{\mu}=0.02,\,\bar{\sigma}=0.2,\,v^h=50,\,v^l=70,\,
\gamma^h=0.2,\,\gamma^l=0.02$, and the terminal condition $g(x)=\min\{x_1,\dots,x_{100}\}$  for $x=( x_1, \dots, x_{ 100 } ) \in \R^{ 100 }$.
Fig.\ \ref{fig:defaultrisk} shows the mean and the standard deviation of $\theta_{u_0}$ 
as an approximation of $ u(t{=}0,x{=}(100,\dots, 100))$, with the final relative error being
$ 0.46\% $.
The not explicitly known ``exact'' solution of \eqref{eq:PDE_defaultrisk} at $t=0,~x=(100,\dots, 100)$ 
has been approximately computed by means of the multilevel Picard method \cite{Weinan2016}:
$u(t{=}0, x{=}(100, \dots, 100)) \approx 57.300$. In comparison, if we do not consider the default risk, 
we get $\tilde{u}(t{=}0,x{=}(100,\dots, 100))\approx60.781$. 
In this case, the model becomes linear and can be solved using straightforward Monte Carlo
methods.
However, neglecting default risks results in a considerable error in the pricing, 
as illustrated above.
The deep BSDE method allows  us
to rigorously incorporate default risks into pricing models.  This in turn
makes it possible to
 evaluate financial derivatives with substantial lower risks for the involved parties 
and the societies.

\subsection*{Hamilton-Jacobi-Bellman (HJB) Equation}
The term ``curse of dimensionality'' was first used explicitly by
Richard Bellman in the context of dynamic programming \cite{Bellman1957}, which has now become
the cornerstone in many areas such as 
economics, behavioral science, computer science, and even biology, where
intelligent decision making is the main issue. 
In the context of game theory where there are multiple players, 
each player has to solve a high-dimensional HJB type equation
in order to find his/her optimal strategy. 
In a dynamic resource allocation problem involving multiple entities with uncertainty,
the dynamic programming principle also leads to a high-dimensional HJB
equation \cite{Powell2011}
for the value function.
Until recently these high-dimensional PDEs have basically remained intractable.
We now demonstrate below that the deep BSDE method is an effective
tool for dealing with these high-dimensional problems.

We consider a classical linear-quadratic-Gaussian (LQG) control problem in 100 dimension:
\begin{equation}
    dX_t = 2\sqrt{\lambda}\,m_t\,dt+\sqrt{2}\,dW_t
\end{equation}
with $ t \in [0,T] $, $ X_0 = x $, and with the cost functional
$ J( \{ m_t \}_{ 0 \leq t \leq T } ) = 
  \E\big[ 
    \int_0^T \|m_t\|^2 \, dt + g(X_T)
  \big]
$.
Here $ \{ X_t \}_{ t \in [0,T] } $ is the state process, 
$ \{ m_t \}_{ t \in [0,T] } $ 
is the control process, 
$ \lambda $ is a positive constant representing the ``strength'' of the control, and $ \{ W_t \}_{ t \in [0,T] } $ is a standard Brownian motion. 
Our goal is to minimize the cost functional through the control process. 
The HJB equation for this problem is given by
\begin{equation}
  \label{eq:PDE_HJB}
     \frac{ \partial u}{ \partial t } ( t, x )
  + \Delta u (t,x) - \lambda \|\nabla u(t,x) \|^2 = 0
\end{equation}
(see e.g., Yong \& Zhou~\cite[Chapter 4]{Yong1999}).
The value of the solution $ u(t,x) $ of \eqref{eq:PDE_HJB} at $t=0$ represents the optimal cost when the 
state starts from $x$. 
Applying It\^{o}'s formula, one can show that the exact solution of \eqref{eq:PDE_HJB}  with the terminal condition
$ u(T, x) = g(x) $ admits the explicit formula
\begin{equation}
\label{eq:HJB_formula}
  u(t,x) = - \frac{ 1}{ \lambda } 
  \ln\!\bigg( 
    \E\Big[ 
      \exp\!\Big( 
         - \lambda g( x + \sqrt{ 2 }W_{ T - t }  )
      \Big) 
    \Big] 
  \bigg).
\end{equation}
This can be used to test the accuracy of the proposed algorithm.

\begin{figure}[ht]
\centering
\setcounter{subfigure}{0}
\subfigure{\includegraphics[width = 7.9cm]{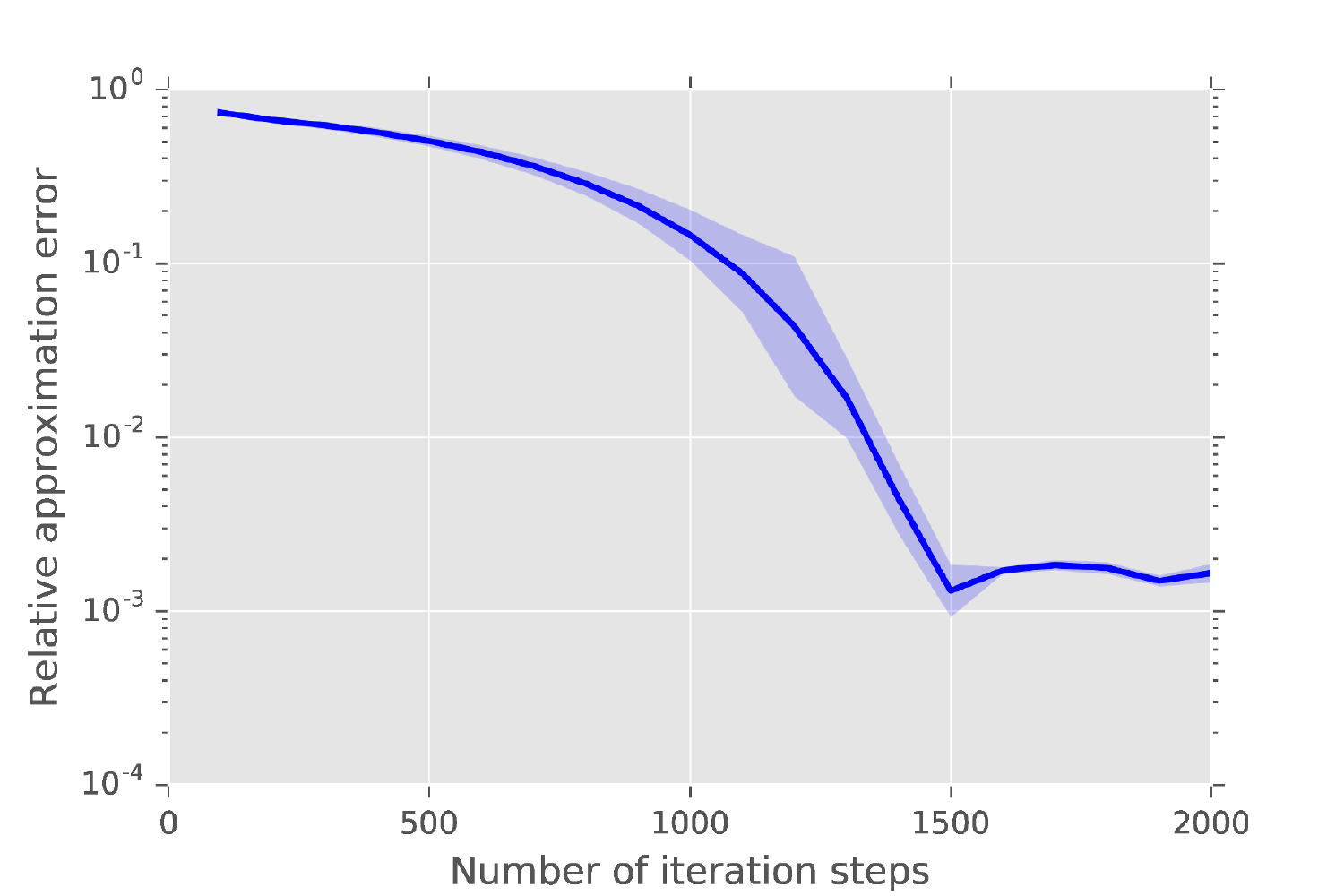}}
\subfigure{\includegraphics[width =7.9cm]{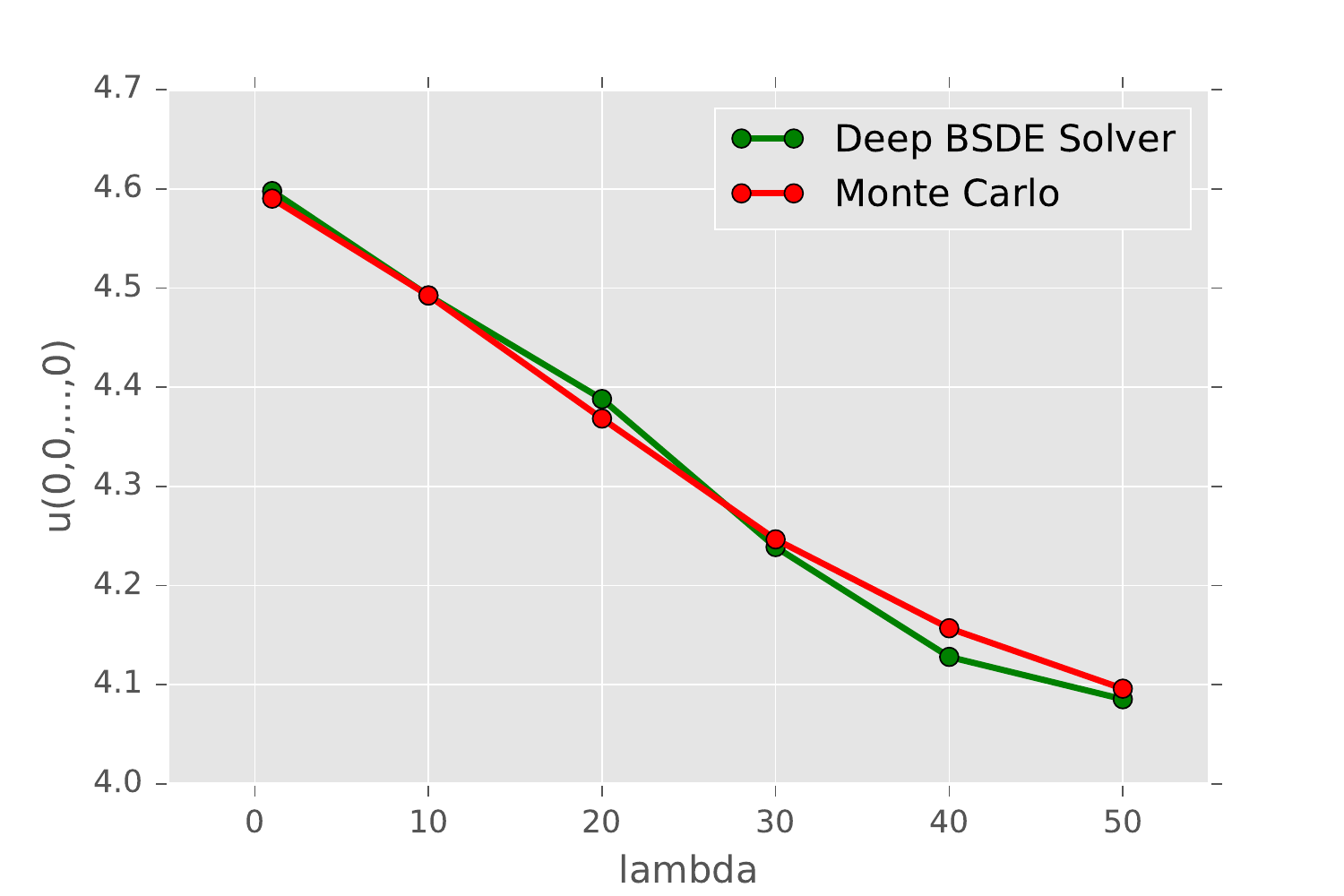}}
\caption{\textbf{Top}: Relative error of the deep BSDE method for $ u( t{=}0, x{=}(0,\dots,0) )$ when $ \lambda = 1 $
against the number of iteration steps in the case of the $ 100 $-dimensional Hamilton-Jacobi-Bellman equation \eqref{eq:PDE_HJB} with $ 20 $ equidistant time steps ($ N{=}20 $) and learning rate $ 0.01 $.
The shaded area depicts the mean $\pm$ the standard deviation of the relative error for 5 different runs.
The deep BSDE method achieves a relative error of size $ 0.17\% $ in a runtime of $ 330 $ seconds. 
\textbf{Bottom}: Optimal cost $u(t{=}0,x{=}(0,\dots,0))$ against different values of $\lambda$ in the case of the $ 100 $-dimensional Hamilton-Jacobi-Bellman equation \eqref{eq:PDE_HJB}, obtained by the deep BSDE method and classical Monte Carlo simulations of \eqref{eq:HJB_formula}.
\label{fig:HJB}}
\end{figure}

We solve the Hamilton-Jacobi-Bellman equation \eqref{eq:PDE_HJB} in the $100$-dimensional case with $g(x) = \ln\left((1 + \| x \|^2)/2\right)$ for $ x \in \R^{ 100 } $.
Fig.\ \ref{fig:HJB} (a) shows the mean and the standard deviation of the relative 
error for $ u(t{=}0,x{=}(0,\dots, 0))$ in the case $\lambda=1$.
The deep BSDE method achieves a relative error of $ 0.17\% $ in a runtime of $ 330 $ seconds
on a Macbook Pro.
We also use the BSDE method to approximatively calculate 
the optimal cost $u(t{=}0,x{=}(0,\dots,0))$ against different 
values of $ \lambda $; see Fig.\ \ref{fig:HJB} (b). 
The curve in Fig.\ \ref{fig:HJB} (b)
clearly confirms the intuition that the optimal cost decreases as 
the control strength increases. 

\subsection*{Allen-Cahn Equation}
\begin{figure}[ht]
\centering
\setcounter{subfigure}{0}
\subfigure{\includegraphics[width=7.9cm]{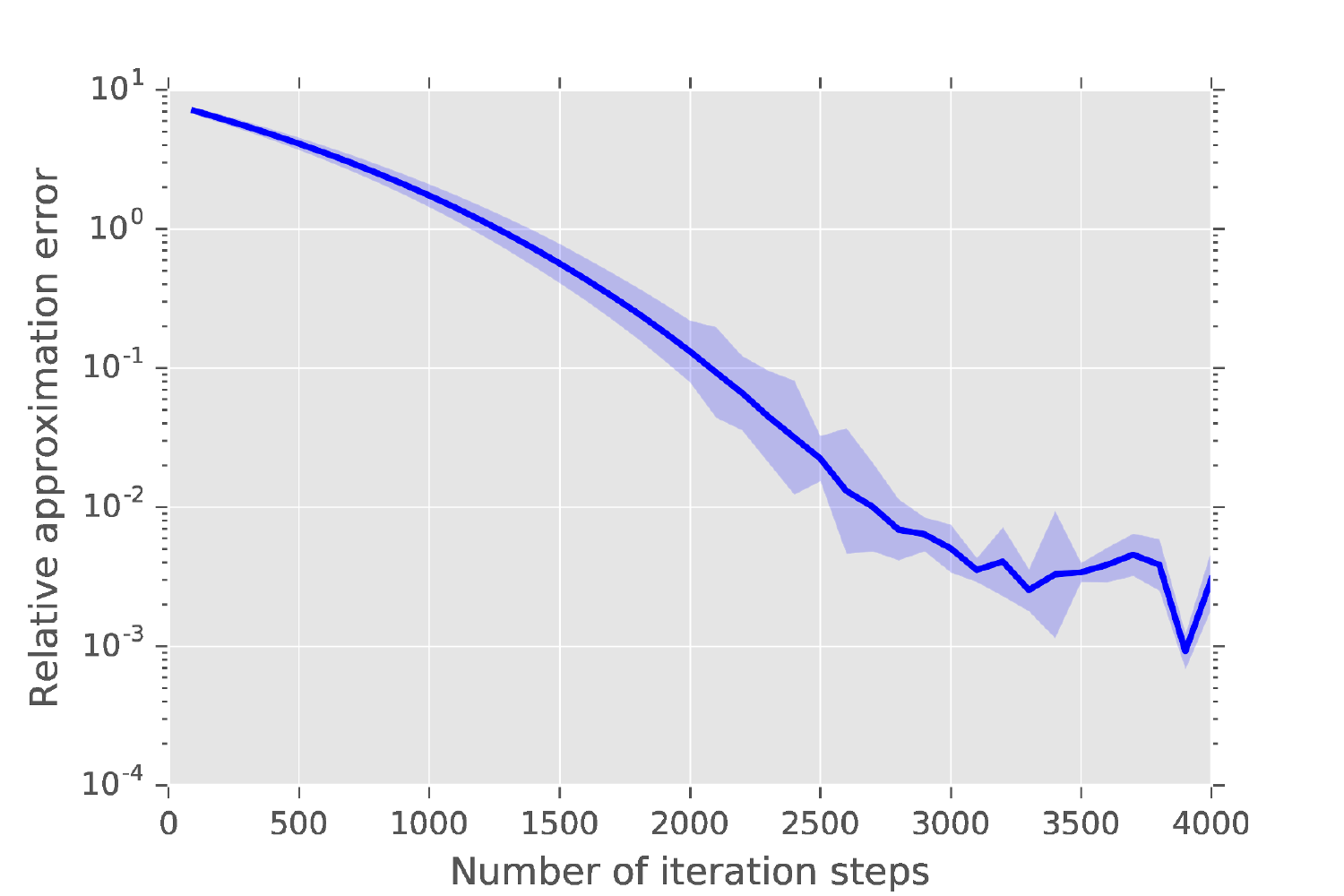}}
\subfigure{\includegraphics[width=7.9cm]{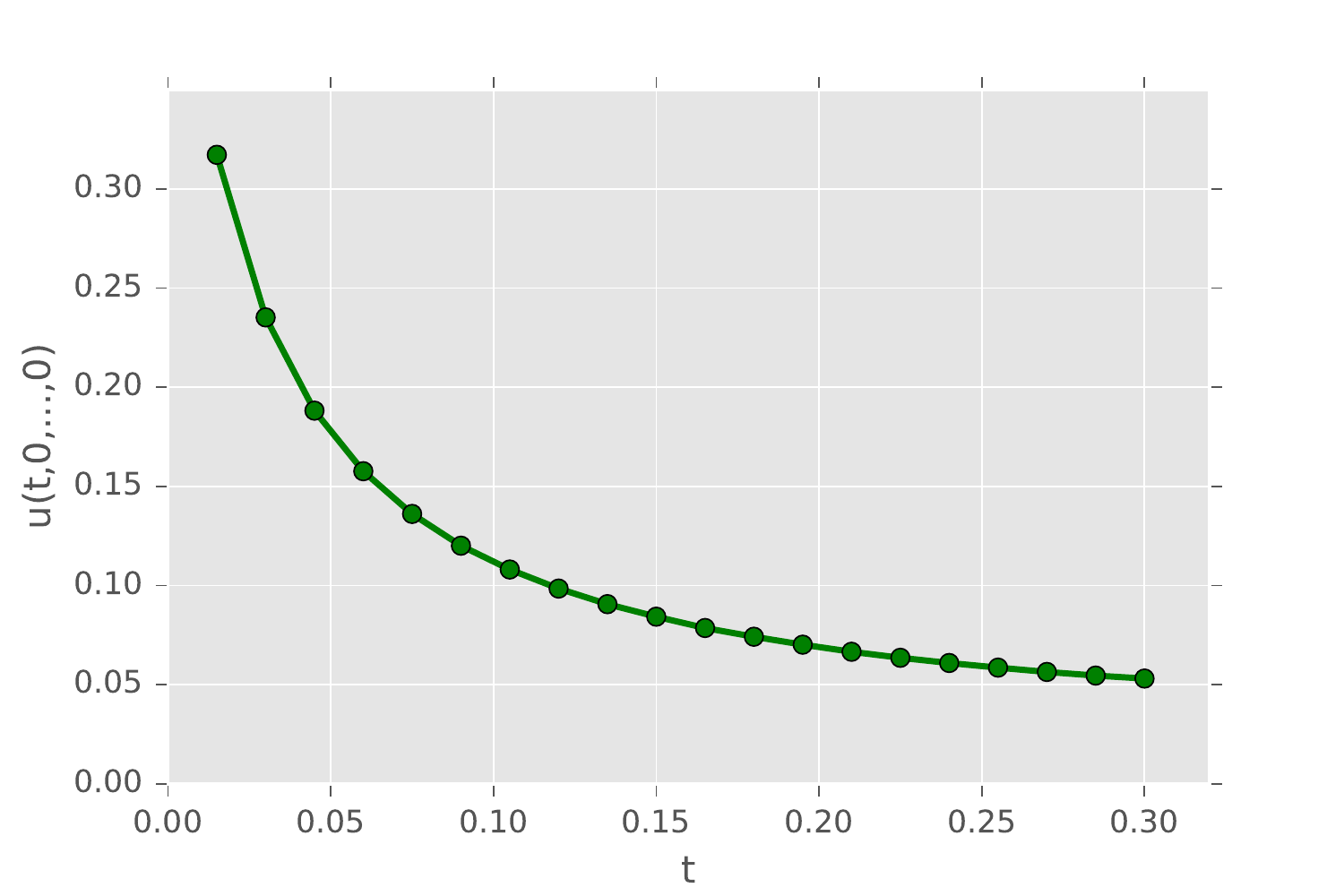}}
\caption{\textbf{Top}: Relative error of the deep BSDE method for $u(t{=}0.3,x{=}(0,\dots,0))$ against the number of iteration steps in the case of the $ 100 $-dimensional Allen-Cahn equation \eqref{eq:PDE_allencahn} with $ 20 $ equidistant time steps ($ N{=}20 $) and learning rate $ 0.0005 $.
The shaded area depicts the mean $\pm$ the standard deviation of the relative error for $ 5 $ different runs. The deep BSDE method achieves a relative error of size $ 0.30\% $ in a runtime of $ 647 $ seconds. 
\textbf{Bottom}: Time evolution of $u(t,x{=}(0,\dots,0))$ for $t\in[0,0.3]$ in the case of the $ 100 $-dimensional Allen-Cahn equation \eqref{eq:PDE_allencahn} computed
by means of the deep BSDE method.
\label{fig:allen_cahn}}
\end{figure}

The Allen-Cahn equation is a reaction-diffusion equation that arises in physics, serving as a 
prototype for the modeling of phase separation and order-disorder transition (see, e.g., \cite{Emmerich2003}).
Here we consider a typical Allen-Cahn equation with the ``double-well potential'' in 100-dimensional space
\begin{equation}
  \label{eq:PDE_allencahn}
     \frac{ \partial u}{ \partial t } ( t, x )
  =  \Delta u (t,x) + u(t,x) - \left[ u(t,x) \right]^3,
\end{equation}
with the initial condition $u(0,x) = g(x)$,
where $g(x)=1/\left(2 + 0.4\, \| x \|^2\right)$ for $ x \in \R^{ 100 } $. By applying a transformation of the time 
variable $t\mapsto T-t \,\,(T>0)$, we can turn \eqref{eq:PDE_allencahn} into the form of \eqref{eq:PDE} 
such that the deep BSDE method can be used. 
Fig.\ \ref{fig:allen_cahn} (a) shows the mean and the standard deviation of the relative 
error of $ u(t{=}0.3,x{=}(0,\dots, 0))$.
The not explicitly known ``exact'' 
solution of \eqref{eq:PDE_allencahn} at $ t = 0.3 $, $ x =(0,\dots,0) $ 
has been approximatively computed by means of 
the branching diffusion method (see, e.g., \cite{Henry-Labordere2012,Henry-Labordere2014}):
$ u( t{=}0.3, x{=}(0, \dots, 0) ) \approx 0.0528$.
For this $ 100 $-dimensional example PDE, the deep BSDE method achieves a relative error of $ 0.30\% $ 
in a runtime of $ 647 $ seconds on a Macbook Pro.
We also use the deep BSDE method to approximatively compute 
the time evolution of $ u(t,x{=}(0,\dots,0)) $ for $ t \in [0,0.3] $; see Fig.\ \ref{fig:allen_cahn} (b).

\section*{\label{sec4}Conclusions}

The algorithm proposed in this paper opens up a host of new possibilities in several
different areas. For example in economics one can consider many different interacting
agents at the same time, instead of using the ``representative agent'' model.
Similarly in finance, one can consider all the participating instruments at the same time,
instead of relying on ad hoc assumptions about their relationships.
In operational research, one can handle 
the cases with hundreds and thousands of participating entities directly, without the need
to make ad hoc approximations.

It should be noted that although the methodology
presented here is fairly general, we are so far not able to deal with the
quantum many-body problem due to the difficulty in dealing with the Pauli
exclusion principle.


\section*{Materials and Methods}
\subsection*{BSDE Reformulation}
The link between (nonlinear) parabolic PDEs and backward stochastic differential equations (BSDEs) has been extensively investigated in the literature 
(see, e.g., \cite{Pardoux1992, Pardoux1999, ElKaroui1997,Gobet2016}). In particular, Markovian BSDEs give a nonlinear Feynman-Kac representation of some nonlinear parabolic PDEs.
Let
$ ( \Omega, \mathcal{F}, \P ) $ 
be a probability space,
$ W \colon [0,T] \times \Omega \to \R^d $
be a $ d $-dimensional standard
Brownian motion,
$ \{ \mathcal{F}_t \}_{ t \in [0,T] } $ 
be the normal filtration generated by $ \{ W_t \}_{ t \in [0,T] } $.
Consider the following BSDE

{\small{
\begin{empheq}[left=\empheqlbrace]{align}
  \label{eq:FSDE_supp}
  &X_t = \xi + \int_0^t\mu( s, X_s )\, ds +\int_0^t \sigma( s, X_s ) \, dW_s, \\
  \label{eq:BSDE_supp}
  &Y_t = g( X_T ) + \int_t^T f( s, X_s, Y_s, Z_s ) \, ds 
  - \int_t^T ( Z_s )^{ \operatorname{T} } \, dW_s,
\end{empheq}
}}for which we are seeking for a $ \{ \mathcal{F}_t \}_{ t \in [0,T] } $-adapted solution process $ \{ ( X_t, Y_t, Z_t ) \}_{t\in[0,T]} $ 
with values in $ \R^d \times \R \times \R^d $.
Under suitable regularity assumptions on the coefficient functions $ \mu $, $ \sigma $, and $ f $, 
one can prove existence and up-to-indistinguishability uniqueness of solutions (cf., e.g., \cite{Pardoux1992,ElKaroui1997}). 
Furthermore, we have that the nonlinear parabolic partial differential equation \eqref{eq:PDE}
is related to the BSDE \eqref{eq:FSDE_supp}--\eqref{eq:BSDE_supp}
in the sense that
for all $ t \in [0,T] $ 
it holds $ \P $-a.s.\ that
\begin{equation}
\label{eq:nonlinear_Feynman_Kac}
  Y_t = u( t, X_t )
\qquad  
  \text{and}
\qquad 
  Z_t = \sigma^{ \operatorname{T} }( t, X_t ) \, 
  \nabla u( t, X_t ),
\end{equation}
(cf., e.g., \cite{Pardoux1992,Pardoux1999}). 
Therefore, we can compute the quantity $u(0, X_0)$ associated to \eqref{eq:PDE} through $Y_0$ by solving the BSDE \eqref{eq:FSDE_supp}--\eqref{eq:BSDE_supp}.
More specifically, we plug the identities in \eqref{eq:nonlinear_Feynman_Kac} into \eqref{eq:BSDE_supp} and rewrite the equation forwardly to obtain the formula in \eqref{eq:BSDE_explicit}.

Then we discretize the equation temporally and use neural networks to approximate the spacial gradients and finally the unknown function, as introduced in the Section Methodology of the paper.

\subsection*{Neural Network Architecture}
\begin{figure*}[ht]
  \centering
  \includegraphics[width=14cm, height=6cm]{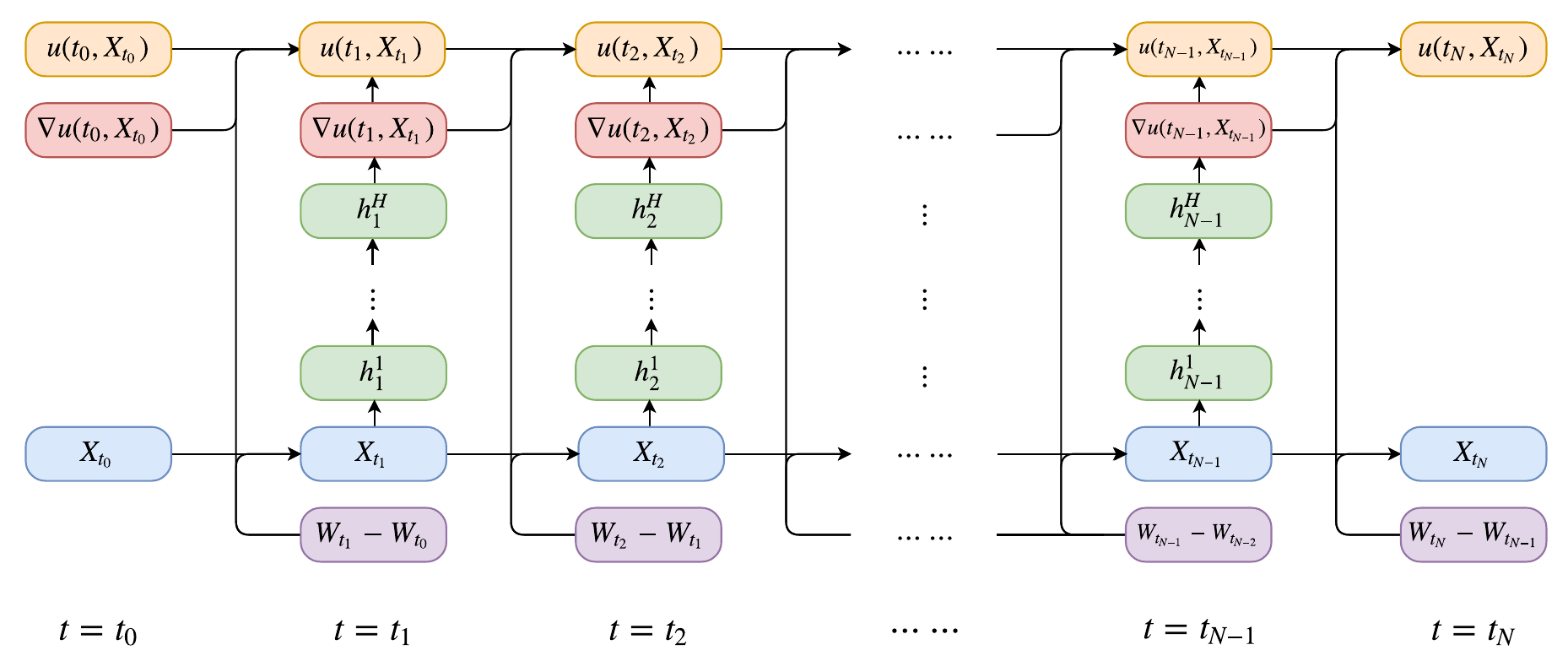}
  \caption{Illustration of the
 network architecture for solving semilinear parabolic PDEs with $H$ hidden layers for 
each sub-network and $N$ time intervals. The whole network has $(H+1)(N-1)$ layers in total that involve free parameters to be optimized simultaneously. Each column for $ t = t_1, t_2,\dots,t_{N-1}$ corresponds to a sub-network at time $t$.
$h^1_n, \dots, h^H_n$ are the intermediate neurons in the sub-network at time $t=t_n$ for $ n= 1, 2, \dots, N - 1$.}
  \label{fig:nn}
\end{figure*}

In this subsection we briefly illustrate the architecture 
of the deep BSDE method.
To simplify the presentation we restrict ourselves in these illustrations 
to the case where the diffusion coefficient $ \sigma $ 
in \eqref{eq:PDE} satisfies that
$ \forall \, x \in \R^d \colon \sigma(x) = \operatorname{Id}_{ \R^d } $.
Fig. \ref{fig:nn} illustrates the network architecture for the deep BSDE method. Note that $\nabla u(t_n,X_{t_n})$ denotes the variable we approximate directly by sub-networks and $ u(t_n,X_{t_n})$ denotes the variable we compute iteratively in the network. There are three types of connections in this network:
\begin{enumerate}
  \item $X_{t_n}\to h_n^1 \rightarrow h_n^2 \rightarrow \cdots \rightarrow h_n^H \rightarrow \nabla u(t_n, X_{t_n})$ is the multilayer feedforward neural network approximating the spatial gradients at time $ t = t_n $. 
  The weights $\theta_n$ of this sub-network are the parameters we aim to optimize.

  \item $(u(t_n,X_{t_n}),\nabla u(t_n, X_{t_n}),W_{t_{n+1}}-W_{t_{n}})\rightarrow 
  u(t_{n+1},X_{t_{n+1}})$ is the forward iteration giving the final output of the network as an approximation of $u(t_N,X_{t_N})$, completely characterized by \eqref{eq:approx_ut}--\eqref{eq:time_discret}. There are no parameters to be optimized in this type of connection.

  \item $(X_{t_n},W_{t_{n+1}}-W_{t_{n}})\rightarrow X_{t_{n+1}}$ is the shortcut connecting blocks at different time, which is characterized by
  \eqref{eq:approx_Xt} and \eqref{eq:time_discret}.
   There are also no parameters to be optimized in this type of connection.
\end{enumerate}
If we use $H$ hidden layers in each sub-network, 
as illustrated in Fig. \ref{fig:nn}, then the whole network has 
$(H+1)(N-1)$ layers in total that involve free parameters to be optimized simultaneously. 

It should be pointed out that the proposed deep BSDE method
can also be employed if we are interested in values 
of the PDE solution $ u $ in a region $ D \subset \R^d $ 
at time $ t = 0 $ instead of at a single space-point $ \xi \in \R^d $. 
In this case we choose $ X_0 = \xi $ to be a non-degenerate 
$ D $-valued random variable and we employ two 
additional neural networks 
parameterized by 
$\{ \theta_{ u_0 }, \theta_{ \nabla u_0} \} $
for approximating the functions
$D \ni x \mapsto u(0,x) \in \R $ and 
$D \ni x \mapsto \nabla u(0,x) \in \R^d$.
Upper and lower bounds for approximation errors of 
stochastic approximation algorithms for PDEs and BSDEs,
respectively, can be found in e.g., \cite{Heinrich2006,Geiss2014,Gobet2016} and the references therein.

\subsection*{Implementation}
We describe in detail the implementation for the numerical examples presented in the paper. 
Each sub-network is fully connected and consists of $ 4 $ layers (except the example in the next subsection), with $ 1 $ input layer ($ d $-dimensional), $ 2 $ hidden layers (both $ d+10 $-dimensional),
and $ 1 $ output layer ($ d $-dimensional).
We choose the rectifier function (ReLU) 
as our activation function.
We also adopted the technique of batch normalization \cite{Ioffe2015} in the sub-networks, right after each linear transformation and before activation. 
This technique accelerates the training by allowing a larger step size and easier parameter initialization. All the parameters are initialized through a normal or a uniform distribution without any pre-training.

We use TensorFlow \cite{Abadi2016} to implement our algorithm with the Adam optimizer \cite{Kingma2015} to optimize parameters. 
Adam is a variant of the SGD algorithm, based on adaptive estimates of lower-order moments. We set the default values for corresponding hyper-parameters as recommended in \cite{Kingma2015} 
and choose the batch size as 64.
In each of the presented numerical examples the means and the standard deviations of the relative $ L^1 $-approximation errors are computed approximatively by means of 5 independent runs of the algorithm 
with different random seeds.
All the numerical examples reported are run on a Macbook Pro with a 2.9GHz Intel Core i5 processor and 16 GB memory.

\subsection*{Effect of Number of Hidden Layers}
The accuracy of the deep BSDE method certainly depends on the number of hidden layers in the sub-network approximation \eqref{eq:approx_grad}. To test this effect, we solve a reaction-diffusion type PDE with different number of hidden layers in the sub-network. The PDE is a high-dimensional version ($d=100$) of the example analyzed numerically in Gobet \& Turkedjiev~\cite{Gobet2017} ($d=2$):
\begin{equation}  
\label{eq:PDE_Gobet}
  \frac{ \partial u}{ \partial t } ( t, x ) + \frac{ 1 }{ 2 } \, \Delta u (t,x) + 
  \min\!\Big\{1 ,\big(u(t,x) - u^*(t,x)\big)^2\Big\}
   = 0,
\end{equation}
in which $u^*(t,x)$ is the explicit oscillating solution
\begin{equation}
  u^*(t,x) = \kappa + \sin\!\big(\textstyle \lambda \sum_{ i = 1 }^d x_i \big) 
      \exp\!\big( \frac{ \lambda^2 d ( t - T ) }{ 2 } \big).
\end{equation}
Parameters are chosen in the same way as in \cite{Gobet2017}: $\kappa=1.6,\,\lambda=0.1,\,T=1$. A residual structure with skip connection is used in each sub-network with each hidden layer having $d$ neurons. We increase the number of hidden layers in each sub-network from $0$ to $4$ and report the relative error in Table \ref{tab:layer}. It is evident that the approximation accuracy increases as the number of hidden layers in the sub-network increases. 

\begin{table}[pht]
\begin{center}
\caption{The mean and standard deviation (std.) of the relative error for the PDE in \eqref{eq:PDE_Gobet}, obtained by the deep BSDE method$^\dagger$ with different number of hidden layers. }
\vspace{8pt}
\label{tab:layer}
\begin{tabular}{@{}c|ccccc}
\toprule
number of layers$^{\dagger\dagger}$ & 29 & 58 & 87 & 116 & 145 
\\  \midrule
mean of relative error & 2.29\% & 0.90\%  & 0.60\% & 0.56\% & 0.53\%\\  
std. of relative error & 0.0026  & 0.0016 & 0.0017 & 0.0017 & 0.0014\\  
\bottomrule
\end{tabular}
\end{center}
\vspace{8pt}
\small{$^\dagger$The PDE is solved until convergence with $30$ equidistant time steps ($N{=}30$) and $40000$ iteration steps. Learning rate is $0.01$ for the first half of iterations and $0.001$ for the second half.\\
$^{\dagger\dagger}$ We only count the layers that have free parameters to be optimized.}
\end{table}

\section*{Acknowledgement}
The work of Han and E is supported in part by Major Program of NNSFC under grant 91130005, DOE grant {DE}-SC0009248 and ONR grant N00014-13-1-0338.

\bibliography{PDE_DL_arxiv}{}
\bibliographystyle{naturemag}

\end{document}